\documentclass{article}

\usepackage{amsfonts}
\usepackage[frenchb,english]{babel}
\usepackage{amsmath}
\usepackage[T1]{fontenc} 
\usepackage{a4wide}
\usepackage{color} 
\usepackage[colorlinks=true, urlcolor=blue,linkcolor=blue, citecolor=blue]{hyperref}

\title{
Asymptotics for random Young diagrams when the word length and alphabet size simultaneously grow to infinity
}
\author{
Jean-Christophe Breton\footnote{Laboratoire de Math\'ematiques, Image et Applications, 
Avenue Michel Cr\'epeau, Universit\'e de La Rochelle, 17042 La Rochelle cedex, France. 
E-mail: jcbreton@univ-lr.fr}
\footnote{Corresponding author}       
\quad and\quad 
Christian Houdr\'e\footnote{
School of Mathematics, 686 Cherry Street,  Georgia Institute of Technology, Atlanta, GA 30332-0160, USA. 
E-mail: houdre@math.gatech.edu}
\footnote{Research supported in part by NSA grant H98230-09-1-0017}
}

\newtheorem{theo}{Theorem}
\newtheorem{prop}[theo]{Proposition}
\newtheorem{corol}[theo]{Corollary}
\newtheorem{lemme}[theo]{Lemma}
\newtheorem{Rem}[theo]{Remark}

\newcommand{\CQFD}{\hfill$\square$}
 
\newcommand{\ind}{\mathbf{1}}

\def\rit{\mathbb{R}}

\def\nit{\mathbb{N}}

\def\pit{\mathbb{P}}
\def\pp{\mathbb{P}}
\def\ee{\mathbb{E}}
\def\E{\mathop{\hbox{\rm I\kern-0.20em E}}\nolimits}
\def\Var{\mathop{\hbox{\rm Var}}\nolimits}

\allowdisplaybreaks 

\begin{document}

\selectlanguage{english}

\maketitle


\begin{abstract}
Given a random word of size $n$ whose letters are drawn independently from an ordered alphabet of size~$m$,
the fluctuations of the shape of the random RSK Young tableaux are investigated, when both $n$ and $m$ converge together to infinity. 
If $m$ does not grow too fast and if the draws are uniform, the limiting shape is the same 
as the limiting spectrum of the GUE.
In the non-uniform case, a control of both highest probabilities will ensure the convergence of the first row of the tableau towards the Tracy-Widom distribution. 
\\
{\bf AMS 2000 Subject Classication.}
Primary: 60F05.
Secondary: 60B12, 60C05, 60F15.
\end{abstract}


\section{Introduction and results}

Let ${\mathcal A}_m=\{\alpha_1<\alpha_2<\cdots<\alpha_m\}$ be an ordered alphabet of size $m$,
and let a word be made of the random letters $X_1^m, \cdots, X_n^m$ (independently) drawn from ${\mathcal A}_m$. 
Recall that the Robinson-Schensted-Knuth (RSK) correspondence associates to a (random) word a pair of (random) Young tableaux of the same shape, having at most $m$ rows ({\it e.g.}, see \cite{Fu} or \cite{Stanley}). 
It is then well known that the length, $V_1(n,m)$, of the top row of these tableaux coincides with the length of the longest (weakly) increasing subsequence of $X_1^m, \dots, X_n^m$. 
The behavior of $V_1(n,m)$ when $n$ and/or $m$ go to $+\infty$ 
and its connections to various mathematical subfields ({\it e.g.} random matrices, queueing theory, percolation theory)
have been investigated in numerous papers (\cite{BDJ}, \cite{BS}, \cite{BM}, \cite{GW}, \cite{ITW1}, \cite{ITW2}, \cite{Johansson}, \cite{TW01}, $\dots$). 
For instance, appropriately renormalized and for uniform draws, $V_1(n,m)$ converges in law, as $n$ goes to infinity and $m$ is fixed, to the largest eigenvalue of a $m\times m$ matrix from the traceless Gaussian Unitary Ensemble (GUE). 
More generally (see \cite{Johansson}),  when $n\to+\infty$ (and $m$ is fixed), 
the shape of the whole Young tableaux associated to a uniform random word converges, after renormalization,
 to the law of the spectrum of a $m\times m$ traceless GUE matrix. 
For different random words such as non uniform or Markovian ones, 
things are more involved (\cite{ITW1}, \cite{ITW2}, \cite{HL2}, \cite{HX}, \cite{CG}). 

For independently and uniformly drawn random words, the following result holds
where, below and in the sequel, $\Rightarrow$ stands for convergence in distribution.  
\begin{theo}
\label{theo:HL2}
Let $V_k(n,m)=\sum_{i=1}^k R_n^i$ be the sum of the lengths $R_n^i$ of the first $k$ rows of the Young tableau. 
Then, 
\begin{equation}
\label{eq:HL2}
\left(\frac{V_k(n,m)-kn/m}{\sqrt n}\right)_{1\leq k\leq m} 
\hskip -.5cm\Rightarrow 
\frac{\sqrt{m-1}}m \left(\max_{{\bf t}\in I_{k,m}} \sum_{j=1}^k\sum_{l=j}^{m-k+j}
\big(\hat B^{l}(t_{j,l})-\hat B^{l}(t_{j,{l-1}})\big)\right)_{1\leq k\leq m}
\hskip -1cm,
\end{equation}
where $(\hat B^1, \dots, \hat B^m)$ is a multidimensional Brownian motion with covariance matrix having dia\-gonal 
terms equal to 1 and off-diagonal terms equal to $-1/(m-1)$,
and where $I_{k,m}$ is defined by 
\begin{align*}
I_{k,m}=\big\{{\bf t}=(t_{j,l}\: :\: 1\leq j\leq k, 0\leq l\leq m) \::\: t_{j,j-1}=0, t_{j, m-k+j}=1, 1\leq j\leq k,&\\
t_{j,l-1}\leq t_{j,l}, 1\leq j\leq k, 1\leq l\leq m-1; \ 
t_{j,l}\leq t_{j-1,l-1
}, 2\leq j\leq k, 2\leq l\leq m
\big\}.&
\end{align*}
\end{theo}
Here, and in the sequel, the rows beyond the height of the tableau are considered to be of length zero. 
Letting ${\bf \Theta}_k:\rit^k\to\rit^k$ be defined via $({\bf \Theta}_k({\bf x}))_j=\sum_{i=1}^j x_i$, $1\leq 
j\leq k$,
then the shape of the Young tableau is given by ${\bf \Theta}_m^{-1}((V_1(n,m), \dots, V_m(n,m))^t)=(R_n^1,\dots, R_n^m)^t$.
Moreover, let $\big(\lambda_{{GUE},m}^{1,0},\lambda_{{ GUE},m}^{2,0}, \dots, \lambda_{{ GUE},m}^{m,0}\big)$ 
be the spectrum, written in non-increasing order, of a $m\times m$ traceless element of the GUE, 
when the GUE is equipped with the measure
$$
\frac 1{C_m}\prod_{1\leq i<j\leq m}(x_i-x_j)^2\prod_{j=1}^me^{-x_j^2/2}
$$
and $C_m=(2\pi)^{m/2}\prod_{j=1}^m j!$ (see \cite{Me}). 
An important fact (see \cite{Ba}, \cite{BJ}, \cite{Doumerc}, \cite{GTW}, \cite{HL2}, \cite{OCY}) 
asserts that 
\begin{align}
\nonumber
&\frac{\sqrt{m-1}}{\sqrt m}{\bf \Theta}_m^{-1}\Big(\Big( \max_{{\bf t}\in I_{k,m}}\sum_{j=1}^k\sum_{l=j}^{m-k+j}
\big(\hat B^{l}(t_{j,l})-\hat B^{l}(t_{j,{l-1}})\big)\Big)_{1\leq k\leq m}\Big)\\
\label{eq:spectrum}
&\hskip 2cm \stackrel{{\cal L}}{=}
\big(\lambda_{{GUE},m}^{1,0},\lambda_{{ GUE},m}^{2,0}, \dots, \lambda_{{ GUE},m}^{m,0}\big).
\end{align}
In fact if
$\big(\lambda_{{ GUE},m}^1,\lambda_{{ GUE},m}^2, \dots, \lambda_{{ GUE},m}^m\big)$
is the (ordered) spectrum of a $m\times m$ element of the GUE, then
\begin{equation}
\label{eq:GUE0}
\big(\lambda_{{ GUE},m}^1,\lambda_{{GUE},m}^2, \dots, \lambda_{{ GUE},m}^m\big)
\stackrel{{\cal L}}{=}
\big(\lambda_{{ GUE},m}^{1,0},\lambda_{{ GUE},m}^{2,0}, \dots, \lambda_{{ GUE},m}^{m,0}\big)
+Z_m e_m,
\end{equation}
where $Z_m$ is a centered Gaussian random variable with variance $1/m$, 
independent of the vector $\big(\lambda_{{ GUE},m}^{1,0},\lambda_{{ GUE},m}^{2,0}, \dots, \lambda_{{ GUE},m}^{m,0}\big)$
and where $e_m=(1,1, \dots, 1)$, see \cite{HX} for simple proofs of \eqref{eq:spectrum} and \eqref{eq:GUE0}. 

Finally, recall that, as $m\to+\infty$, the asymptotic behavior of the spectrum of the GUE has been obtained by Tracy and Widom 
(see \cite{TW94},  \cite{TW98} and also Theorem 1.4 in \cite{Johansson}, with slight change in the notation):
\begin{theo}
\label{theo:Johansson}
For each $r\geq 1$, there is a distribution ${\bf F}_{r}$ on $\rit^r$ such that: 
\begin{equation}
\label{eq:Johansson}
\left(m^{1/6}\big(\lambda_{{ GUE},m}^k-2\sqrt{m}\big)\right)_{1\leq k\leq r} 
\Rightarrow {\bf F}_{r}, \quad m\to+\infty.
\end{equation}
\end{theo}
 
\begin{Rem}
\label{rem:Johansson-TW}
{\rm The distribution ${\bf F}_{r}$ is explicitly known (see (3.48) in \cite{Johansson}) 
and its first marginal coincides with the Tracy-Widom distribution.  
}\end{Rem}

Since $Z_m m^{1/6}\Rightarrow 0$ as $m\to+\infty$, taking successively the limits in $n$ and then in $m$, 
\eqref{eq:HL2}--\eqref{eq:Johansson} entail for each $r\geq 1$:
\begin{equation}
\label{eq:nm2}
\lim_{m\to+\infty}\lim_{n\to+\infty}
\left(\frac{V_k(n,m)-kn/m-2k\sqrt{n}}{\sqrt n}\times  m^{2/3}\right)_{1\leq k\leq r} = {\bf F}_{r} {\bf \Theta}_r^{-1}.
\end{equation}
In fact, since $\lim_{m\to +\infty}m^{1/6}Z_m = 0$, in probability, and in view of \eqref{eq:GUE0}, throughout, in studying weak asymptotics,  one will be able to replace the correlated Brownian motions of \eqref{eq:spectrum} by uncorrelated (standard) ones.

Following universality argument in percolation models developed by Bodineau and Martin (\cite{BM}), 
we show below that the limits in $n$ and $m$ in \eqref{eq:nm2} can be explicitly taken simultaneously 
when the size $m$ of the alphabet does not grow too fast with respect to $n$. 
Doing so, we are dealing with growing ordered alphabets and at each step, 
the $n$ letters $X_i^m$, $1\leq i\leq n$, are redrawn  
(and not just the $n$th letter as in the case with the model studied in \cite{HIL}). 
In a way, we are thus giving the fluctuations of the shape of the Young tableau of a random word when the alphabets are growing and are reshuffled. 
In the sequel, $m$ will be a function $m(n)$ of $n$. However in order to 
lighten the notation, we shall still write $m$ instead of $m(n)$. 
A main result of this note is:
\begin{theo}
\label{theo:main}
Let $m$ tend to infinity as $n\to+\infty$ in such a way that $m=o(n^{3/10}(\log n)^{-3/5})$. 
Then for each $r\geq 1$,  
$$
\left(\frac{V_k(n,m)-kn/m-2k\sqrt{n}}{n^{1/2}m^{-2/3}}\right)_{1\leq k\leq r} 
\Rightarrow {\bf F}_{r} {\bf \Theta}_r^{-1}, \quad n\to+\infty.
$$
\end{theo}
Remark \ref{rem:condition}, below, briefly discusses the growth conditions on $m$. 
Since, again, the length of the first row of the Young tableau is 
the length $V_1(n,m)$ of the longest increasing subsequence and since the first marginal of ${\bf F}_r$ is the 
Tracy-Widom distribution $F_{TW}$, 
we have: 

\begin{corol}
\label{corol:TW}
Let $m$ tend to infinity as $n\to+\infty$ in such a way that $m=o(n^{3/10}(\log n)^{-3/5})$. 
Then 
$$
\frac{V_1(n,m)-(n/m)-2n^{1/2}}{n^{1/2}m^{-2/3}} \Rightarrow F_{TW}, \quad n\to+\infty.
$$
\end{corol}

When the independent random letters are no longer uniformly drawn, a similar asymptotic behavior continues to hold for $V_1(n,m)$ as explained next. 
Let the $X_i^m$, $1\leq i\leq n$, be independently  and identically distributed with $\pit(X_1^m=\alpha_j)=p_j^m$, 
let $p_{max}^m=\max_{1\leq j\leq m} p_j^m$, and let also $J(m)=\{j\: :\: p_j^m=p_{max}^m\}=\{j_1, \dots, j_{k(m)}\}$ 
with $k(m)=\mbox{card }(J(m))$. 
Now, from \cite{HL1} and as $n\to+\infty$, the behavior of the first row of the Young tableau in this non-uniform setting is given by: 
\begin{equation}
\label{eq:nonunif1}
\frac{V_1(n,m)-p_{max}^mn}{\sqrt{p_{max}^mn}}\Rightarrow
\frac{\sqrt{1-k(m)p_{max}^m}-1}{k(m)}\sum_{j=1}^{k(m)} B^j(1)
+\max_{\begin{subarray}{l} 0=t_0\leq t_1\leq\dots\\ \leq t_{k(m)-1}\leq t_{k(m)}=1\end{subarray}}
\sum_{l=1}^{k(m)}(B^l(t_l)-B^l(t_{l-1})),
\end{equation}
where $(B^1, \dots, B^{k(m)})$ is a standard $k(m)$-dimensional Brownian motion. 
For the limiting behavior in $m$ of the right-hand side of \eqref{eq:nonunif1}, and as explained next two cases can arise, depending on the number of most probable letters in ${\cal A}_m$.
Setting,
$$
Z_k= \frac 1{k}\sum_{j=1}^{k} B^j(1)
\quad \mbox{ and } \quad
D_{k}=\max_{\begin{subarray}{l} 0=t_0\leq t_1\leq\dots\\ \leq t_{k-1}\leq t_{k}=1\end{subarray}}
\sum_{l=1}^{k}\big(B^l(t_l)-B^l(t_{l-1}))\big), 
$$
and combining \eqref{eq:spectrum}, \eqref{eq:GUE0} and \eqref{eq:Johansson} as well as Remark \ref{rem:Johansson-TW}, when $k=1$,
and since clearly $Z_k\sim {\cal N}(0, 1/k)$, we have: 
\begin{equation}
\label{eq:Dk}
k^{1/6}\big(D_k-2\sqrt{k}\big)\Rightarrow F_{TW}, \quad k\to+\infty. 
\end{equation}

First, let $k(m)$ be bounded. 
Eventually extracting a subsequence, we can assume that $k(m)$ is equal to a fixed $k\in \nit\setminus\{0\}$ 
and since $p_{max}^m\in [0,1]$, we can also assume that $p_{max}^m\to p_{max}$. 
In this case, taking the limit first in $n$ and next in $m$ yields:   
\begin{equation}
\label{eq:lim_2}
\frac{V_1(n,m)-p_{max}^mn}{\sqrt{p_{max}^mn}}
\Rightarrow 
\big( \sqrt{1-kp_{max}}-1\big) Z_k+D_k.
\end{equation}
The limiting distribution on the right-hand side of \eqref{eq:lim_2} depends on $k$. 
For instance for $k=1$, we recover a Gaussian distribution while 
for $k>1$ and specific choice of the $p_{max}^m$ for which $\lim_{m\to+\infty} p_{max}^m=0$, 
we recover \eqref{eq:lim_2} without the Gaussian term.
Thus, in general, when $k(m)$ is bounded, there is  no global asymptotics but only  convergence (to different distributions) 
along subsequences. 

Next, let $k(m)\to+\infty$. In this case, in \eqref{eq:nonunif1}, the Gaussian contribution is negligible. 
Indeed, since $(\sqrt{1-{k(m)}p_{max}^m}-1)^2k(m)^{-2/3}\leq (k(m)p_{max}^m)^2 k(m)^{-2/3}\leq k(m)^{-2/3}\to 0$, 
when $m\to+\infty$:
$$
(\sqrt{1-k(m)p_{max}^m}-1) Z_{k(m)} k(m)^{1/6}
\sim {\mathcal N}(0,(\sqrt{1-{k(m)}p_{max}^m}-1)^2k(m)^{-2/3}\big) \Rightarrow 0.
$$
Hence plugging the convergence result \eqref{eq:Dk} into \eqref{eq:nonunif1} leads to 
\begin{equation}
\label{eq:lim_TW}
\frac{V_1(n,m)-p_{max}^mn-2\sqrt{k(m)p_{max}^mn}}{\sqrt{k(m)p_{max}^mn}} k(m)^{2/3}
\Rightarrow F_{TW}
\end{equation}
where the limit is first taken as $n\to+\infty$ and then as $m\to+\infty$. 
In this non-uniform setting, we have the following counterpart to Corollary \ref{corol:TW} with an additional control on 
the second largest probability for the letters of ${\cal A}_m$. 
More precisely, let $p_{2nd}^m=\max(p_j^m<p_{max}^m\: :\: 1\leq j\leq m)$:  
\begin{theo}
\label{theo:nonunif_TW}
Let the size $m$ of the alphabets vary with $n$ 
and assume that $k(m(n))$, the number of most probable letters in ${\cal A}_m$, goes to infinity when $n\to+\infty$, in such a way that 
$k(m(n))^{7/10}/p_{max}^{3/10}=o(n^{3/10}(\log n)^{-3/5})$.
Assume moreover that
\begin{equation}
\label{eq:condition2}
(p_{2nd}^{m(n)})^2\frac{n^{11/10}}{(\log n)^{1/5}}=o(p_{max}^{m(n)}). 
\end{equation}
Then 
\begin{equation}
\label{eq:nonunif_TW}
\frac{V_1(n,m(n))-p_{max}^{m(n)}n-2\sqrt{k(m(n)) p_{max}^{m(n)}n}}{\sqrt{k(m(n))p_{max}^{m(n)}n}}  k(m(n))^{2/3} \Rightarrow F_{TW}.
\end{equation}
\end{theo}
Let us stress again the fact that in the previous result, $m$ is a function of $n$, with the only requirement that $k(m(n))^{7/10}/(p_{max}^{m(n)})^{3/10}=o(n^{3/10}(\log n)^{-3/5})$.  
Note that in the uniform case, $k(m)=m$ and $p_{max}^m=1/m$ and that in general $1/m\leq p_{max}^m\leq 1/k(m)$. 

Let us now put our results in context, relate them to the current literature, and also describe the main steps in the arguments developed below. 

Bodineau and Martin \cite{BM} showed that the fluctuations of the last-passage directed percolation model with Gaussian {\it iid} weights actually extend to {\it iid} weights with finite $(2+r)$-th moment, $r>0$.
Their arguments rely, in part, on a KMT approximation 
which was already used by Glynn and Whitt \cite{GW} in a related queueing model. 

Here, we closely follow \cite{BM} 
and take advantage of the representation \eqref{eq:spectrum} of the spectrum of a matrix in the GUE. 
Using Brownian scaling in those Brownian functionals, we can mix together $n$ and $m$ in the corresponding limit \eqref{eq:Johansson} (see \eqref{eq:scaling2} below). 
Then, exhibiting an expression similar to \eqref{eq:spectrum}, but with {\it dependent} Bernoulli random variables, for the shape of the Young tableau (see \eqref{eq:Vn}), 
we show via a Gaussian approximation that the Bernoulli functionals stay close to the Brownian functionals (see \eqref{eq:bound1}), so as to share the same asymptotics. 

Since we apply a Gaussian approximation to Bernoulli random variables with strong integrability property, 
the strong approximation can be made more precise than in \cite{BM}.
But, this is not enough to obtain the fluctuations for $m$ of larger order. 
Actually the Gaussian approximation is responsible for the condition $m=o\big(n^{3/10}(\log n)^{-3/5}\big)$, 
which comes short of the corresponding polynomial order condition $m=o(n^{3/7})$ obtained in \cite{BM}. 
However in contrast to \cite{BM}, the stronger integrability property of the Bernoulli random variables and the stronger condition on $m$
are required to control the constants appearing in the Gaussian approximation applied to a triangular scheme of different distributions. 

Using Skorohod embedding, Baik and Suidan \cite{BS} derived, independently of \cite{BM}, similar convergence results 
(see Theorem 2 in \cite{BS}), under the condition $m=o\big(n^{3/14}\big)$. 
See also \cite{Su} for related results (under $m=o\big(n^{1/7}\big)$) in percolation models using functional methods in the CLT.

Finally, note that \cite{BM, BS, Su} deal with percolation models with {\it iid} 
random variables under enough polynomial integrability. 
In our setting, the lengths of the rows of the Young tableaux associated to random words 
are expressed in terms of dependent (exchangeable in the uniform case) Bernoulli random variables. 
We are thus working with much more specific random variables but without complete independence. 

The paper is organized as follows: Section \ref{sec:proof_main} is devoted to the proof of Theorem \ref{theo:main}, 
while we sketch the changes needed to prove Theorem \ref{theo:nonunif_TW} in Section \ref{sec:nonunif_TW}. 
We conclude in Section \ref{sec:remarks} with some remarks on the convergence of whole shape of Young tableaux when the draws are non uniform.


\section{Proof of Theorem \ref{theo:main}}
\label{sec:proof_main}

\hskip 15pt
{\bf Brownian scaling.}
Let $(B^{l}(s))_{s\geq 0}$, $1\leq l\leq m$, be independent standard Brownian motions. 
For $s>0$, $m\geq 1$ and $k\geq 1$, let
\begin{equation}
\label{eq:percolation2}
L_k(s,m)=\sup_{{\bf t}\in I_{k,m}(s)} \sum_{j=1}^k\sum_{l=j}^{m-k+j}\big(B^{l}(t_{j,l})-B^{l}(t_{j,l-1})\big),
\end{equation}
where $I_{k,m}(s)=\{s{\bf t}, {\bf t}\in I_{k,m}\}$. 
For $k=1$, $L_1(s,m)$ coincides with the Brownian percolation model used in \cite{BM}, 
see also \cite{GW} for a related queueing model. 
For $s=1$, $
{\bf \Theta}_m^{-1}((L_k(1,m))_{1\leq k\leq m})$ has the same law as the spectrum of a $m\times m$ GUE matrix, see \cite{Doumerc} and \cite{HX}.\\ 
Since $(L_1(\cdot,m), \dots, L_m(\cdot,m))$ is a continuous function of $B^{1}, \dots, B^{m}$, which are independent;
Brownian scaling entails: 
\begin{equation}
\label{eq:scaling1}
\big(L_1(s,m), \dots,L_m(s,m)\big)\stackrel{{\cal L}}{=}\sqrt s \big(L_1(1,m), \dots, L_m(1,m)\big).
\end{equation}
Plugging \eqref{eq:scaling1} into \eqref{eq:Johansson} yields, as $m\to+\infty$, 
\begin{equation}
\label{eq:scaling2}
\left(\frac{L_k(n,m)-2k\sqrt{nm}}{n^{1/2}m^{-1/6}}\right)_{1\leq k\leq r}
\Rightarrow {\bf F}_{r} {\bf \Theta}_r^{-1}. 
\end{equation}

{\bf Combinatorics.}
Let 
$$
X_{i,j}^m=\left\{\begin{array}{ll}
1&\mbox{if } X_i^m=\alpha_j\\
0&\mbox{otherwise,}
\end{array} \right .
$$
be Bernoulli random variables with parameter $\pit(X_i^m=\alpha_j)=1/m$ and variance $\sigma_m^2=(1/m)(1-1/m)$. 
For a fixed $1\leq j\leq m$, the $X_{i,j}^m$s are independent and identically distributed while for $j\not=j'$, 
$(X_{1,j}^m, \dots, X_{n,j}^m)$ and $(X_{1,j'}^m, \dots, X_{n,j'}^m)$ are identically distributed but no longer independent. 

Recall again that the length of the first row of the Young tableau of a random word is the length of the longest (weakly) increasing subsequence of $X_1^m, \dots, X_n^m$. 

Let $S_k^{m,j}=\sum_{i=1}^k X_{i,j}^m$ be the number of occurences of $\alpha_j$ among $(X_i^m)_{1\leq i\leq k}$. 
An increasing subsequence of $(X_i^m)_{1\leq i\leq k}$ consists of successive blocks, each one made of an identical letter, with the sequence of letters representing each block being strictly increasing.
Since for $1\leq k<l\leq n$ the number of occurences of $\alpha_j$ among $(X_i^m)_{k\leq i\leq l}$ is $S_l^{m,j}-S_k^{m,j}$, and it follows that:
\begin{equation}
\label{eq:combinatoire_V_1}
V_1(n,m)
=\max_{\begin{subarray}{l}0=l_0\leq l_1\leq\cdots\\\leq l_{m-1}\leq l_m=n\end{subarray}}
\big[(S_{l_1}^{m,1}-S_0^{m,1})+(S_{l_2}^{m,2}-S_{l_1}^{m,2}) 
+\cdots+(S_{n}^{m,m}-S_{l_{m-1}}^{m,m})\big],
\end{equation}
with the convention that $S_0^{m,1}=0$. More involved combinatorial arguments yield the following expression for $V_k(n,m)$ (see Theorem 5.1 in \cite{HL2}):
\begin{equation}
\label{eq:combinatoire_V_k}
V_k(n,m)=\max_{{\bf k}\in J_{k,m}(n)} \sum_{j=1}^k\sum_{l=j}^{m-k+j}\big(S_{k_{j,l}}^{m,l}-S_{k_{j,l-1}}^{m,l}\big),
\end{equation}
where
\begin{align*}
&J_{r,m}(n)=\big\{{\bf k}=(k_{j,l}\: :\: 1\leq j\leq r, 0\leq l\leq m) \::\: k_{j,j-1}=0, k_{j, m-r+j}=n, 1\leq j\leq r,\\
&\hskip 2cm k_{j,l-1}\leq k_{j,l}, 1\leq j\leq r, 1\leq l\leq m-1;\: k_{j,l}\leq k_{j-1,l-1}, 2\leq j\leq r, 1\leq l\leq m\big\}.
\end{align*}
For ${\bf t}\in I_{r,m}(n)$, set $[{\bf t}]=\big([t_{j,l}]\: : \: 1\leq j\leq n, 0\leq l\leq m\big)\in J_{r,m}(n)$ 
and thus
\begin{equation}
\label{eq:Vn}
V_k(n,m)=\sup_{{\bf t}\in I_{k,m}(n)} \sum_{j=1}^k\sum_{l=j}^{m-k+j}\big(S_{[t_{j,l}]}^{m,l}-S_{[t_{j,l-1}]}^{m,l}\big),
\end{equation}
which is to be compared with \eqref{eq:percolation2} for Brownian functionals. 

{\bf Centering and reducing.}
Let $\widetilde{X}_{i,j}^m=(X_{i,j}^m-1/m)/\sigma_m$ 
and $\widetilde{S}_k^{m,l}=\sum_{i=1}^k \widetilde{X}_{i,l}^{m}$, 
and replacing $X_{i,j}^m$ by $\widetilde{X}_{i,j}^m$, similarly define $\widetilde V_k(n,m)$. 
Clearly, $V_k(n,m)=\sigma_m \widetilde{V}_k(n,m)+kn/m$, hence,

\begin{align*}
&\frac{V_k(n,m)-kn/m-2k\sqrt{n}}{\sqrt n}\times  m^{2/3}\\
&=\frac{\sigma_m\widetilde{V}_k(n,m)-2k\sqrt{n}}{\sqrt n}\times  m^{2/3}\\
&=\frac{\widetilde{V}_k(n,m)-2k\sqrt{n}\sigma_m^{-1}}{\sqrt n}\times (\sigma_m m^{2/3})\\
&=\frac{\widetilde{V}_k(n,m)-2k\sqrt{nm}+2k\sqrt{n}(\sigma_m^{-1}-m^{1/2})}{n^{1/2}m^{-1/6}}
\times (m^{1/2}\sigma_m).
\end{align*}
Note that $\sigma_m^{-1}-m^{1/2}\sim 1/\sqrt m$, and that $m^{1/6}m^{1/2}\sigma_m\sim m^{1/6}$, and so
the limit under study is the same as that of
\begin{equation}
\label{eq:V_k}
\frac{\widetilde{V}_k(n,m)-2k\sqrt{nm}}{n^{1/2}m^{-1/6}}.
\end{equation}

\medskip
{\bf Bound.}
Next and as \cite{BM}, we bound the difference between $\widetilde{V}_k(n,m)$ and $L_k(n,m)$. 
This bound holds true for any Brownian motions $(B_t^{m,j})_{t\geq 0}$ 
but it will only be correctly controlled for a special choice of the Brownian motions 
and for copies of the random variables $\widetilde{X}_{i,j}^m$
given by a coupling (using a strong approximation result, see Proposition~\ref{prop:KMT} below).
\begin{align}
\nonumber
&\left|\widetilde{V}_k(n,m)-L_k(n,m)\right|\\
\nonumber
&=\left|
\sup_{{\bf t}\in I_{k,m}(n)} \sum_{j=1}^k\sum_{l=j}^{m-k+j}\big(\widetilde S_{[t_{j,l}]}^{m,l}-\widetilde S_{[t_{j,l-1}]}^{m,l}\big)
-\sup_{{\bf t}\in I_{k,m}(n)}\sum_{j=1}^k\sum_{l=j}^{m-k+j}\big(B^{l}(t_{j,l})-B^{l}(t_{j,l-1})\big)\right|\\
\nonumber
&\leq\sup_{{\bf t}\in I_{k,m}(n)}
\left|
 \sum_{j=1}^k\sum_{l=j}^{m-k+j}\big(\widetilde S_{[t_{j,l}]}^{m,l}-\widetilde S_{[t_{j,l-1}]}^{m,l}\big)
-\sum_{j=1}^k\sum_{l=j}^{m-k+j}\big(B^{l}(t_{j,l})-B^{l}(t_{j,l-1})\big)\right|\\
\nonumber
&=\sup_{{\bf t}\in I_{k,m}(n)}
\left|
 \sum_{j=1}^k\sum_{l=j}^{m-k+j}\big(\widetilde S_{[t_{j,l}]}^{m,l}-B^{l}(t_{j,l})\big)
-\sum_{j=1}^k\sum_{l=j}^{m-k+j}\big(\widetilde S_{[t_{j,l-1}]}^{m,l}-B^{l}(t_{j,l-1})\big)\right|\\
\nonumber
&=\sup_{{\bf t}\in I_{k,m}(n)}
\Bigg|
 \sum_{j=1}^k\sum_{l=j}^{m-k+j}
 \Big(\big(\widetilde S_{[t_{j,l}]}^{m,l}-B^{l}([t_{j,l}])\big)+\big(B^{l}([t_{j,l}])-B^{l}(t_{j,l})\big)\\
\nonumber &\hskip 3cm
-\big(\widetilde S_{[t_{j,l-1}]}^{m,l}-B^{l}([t_{j,l-1}])\big)-\big(B^{l}([t_{j,l-1}])-B^{l}(t_{j,l-1})\big)\Big)\Bigg|\\
\nonumber
&\leq\sup_{{\bf t}\in I_{k,m}(n)}
\Bigg\{
 \sum_{j=1}^k\sum_{l=j}^{m-k+j}
 \Big(|\widetilde S_{[t_{j,l}]}^{m,l}-B^{l}([t_{j,l}])|+|B^{l}([t_{j,l}])-B^{l}(t_{j,l})|\\
\nonumber &\hskip 3cm
+|\widetilde S_{[t_{j,l-1}]}^{m,l}-B^{l}([t_{j,l-1}])|+|B^{l}([t_{j,l-1}])-B^{l}(t_{j,l-1})|\Big)\Bigg\}\\
\label{eq:bound1}
&\leq 2k\sum_{l=1}^{m}\left(Y_n^{m,l}+W_n^{l}\right),
\end{align}
where we set 
$$
Y_n^{m,l}=\max_{1\leq i\leq n} |\widetilde S_i^{m,l}-B^{l}(i)|
\quad \mbox{ and } \quad 
W_n^{l}=\sup_{\begin{subarray}{c} 0\leq s,t\leq n\\ |s-t|\leq 1\end{subarray}}|B^{l}(s)-B^{l}(t)|.
$$

\medskip
{\bf Gaussian approximation.}
From now on, we assume that for each $n$ and $l\in [1,m]$ (recall that $m=m(n)$), 
the random variables $\widetilde{X}_{i,l}^m$, $1\leq i\leq n$, and the Brownian motion $(B^{l}(s))_{s\in [0,n+1]}$, 
appearing in $Y_n^{m,l}$ and $W_n^{l}$ (rewritten as $(B^{m,l}(s))_{s\in [0,n+1]}$), are given by the following result, 
which is a compilation of strong approximation results of Koml\'os, Major, Tusn\'ady and of Sakhanenko and for which we refer to \cite{L} (Th. 2.1, Cor 3.2) and the references therein. 
In the sequel, we write $B^{m,l}$ and $W_n^{m,l}$, instead of $B^{l}$ and $W_n^{l}$,
to insist on the dependence in $m$ of the random variables given by the forthcoming proposition. 
\begin{prop}
\label{prop:KMT}
Let $(X_n)_{n\geq 1}$ be a sequence of {\it iid} random variables with common distribution $F$ having finite exponential moments. 
Then, on a common probability space and for every $N$, 
one can construct a sequence $(\widetilde{X}_n)_{1\leq n\leq N}$ having the same law as $(X_n)_{1\leq n\leq N}$,
and independent Gaussian variables $(Y_n)_{1\leq n\leq N}$ having same expectations and variances as $(X_n)_{1\leq n\leq N}$ such that for every $x>0$:
$$
\pit\left(\max_{1\leq k\leq N}\left|\sum_{j=1}^k \widetilde{X}_j-\sum_{j=1}^k Y_j\right|\geq x\right)
\leq(1+c_2(F)N^{1/2}) \exp(-c_1(F)x),
$$
where $c_1(F)$ and $c_2(F)$ are positive constants (depending on $F$).
Moreover $c_1(F)=c_3\lambda(F)$ and $c_2(F)=\lambda(F) \Var(X_1)^{1/2}$, 
where $c_3$ is an absolute constant and $\lambda(F)$ is given by
$$
\lambda(F)=\sup\left\{\lambda>0 \: :\: \lambda \ee[|X_1-\ee[X_1]|^3\exp(\lambda|X_1-\ee[X_1]|)]\leq \ee[|X_1-\ee[X_1]|^2]\right\}.
$$
\end{prop}

The strong approximation entails the following bound for the tail of $Y_n^{m,l}$:
\begin{equation}
\label{eq:KMT}
\pit(Y_n^{m,l}\geq x)\leq (1+c_2(m)n^{1/2})\exp(-c_1(m)x),
\end{equation}
where $c_1(m)=c_3\lambda(\widetilde{X}^m_{1,1})$ and 
$c_2(m)=\lambda(\widetilde{X}^m_{1,1}) \Var(\widetilde{X}^m_{1,1})^{1/2}$. 
Observe that  $\lambda(\widetilde{X}^m_{1,1})=\sigma_m\lambda(X^m_{1,1}-\ee[X^m_{1,1}])$ 
and note that $\lambda(X^m_{1,1})\in [2^{-1},2]$. 
Indeed, for $\lambda\geq 2$,
\begin{eqnarray*}
\ee[|X^m_{1,1}-\ee[X^m_{1,1}]|^2]&=&\frac 1m\Big(1-\frac 1m\Big)\\
&\leq& \frac 1m\Big(1-\frac 1m\Big) \frac{\lambda}2\\
&\leq&\frac 1m\Big(1-\frac 1m\Big) \Big(1-\frac 2m+\frac2{m^2}\Big)\lambda\\
&=&\lambda\ee[|\widetilde{X}^m_{1,1}-\ee[X^m_{1,1}]|^3]\\
&\leq&\lambda\ee[|X^m_{1,1}-\ee[X^m_{1,1}]|^3\exp(\lambda|\widetilde{X}^m_{1,1}-\ee[X^m_{1,1}]|)],
\end{eqnarray*}
while, since $|X^m_{1,1}-\ee[X^m_{1,1}]|\leq 1$,
$$
\frac 12\ee\Big[|X^m_{1,1}-\ee[X^m_{1,1}]|^3\exp\Big(\frac 12|X^m_{1,1}-\ee[X^m_{1,1}]|\Big)\Big]
\leq \frac 12\exp\Big(\frac 12\Big)\ee[|X^m_{1,1}-\ee[X^m_{1,1}]|^2]
\leq \ee[|X^m_{1,1}-\ee[X^m_{1,1}]|^2].
$$
Thus, $c_1(m)$ and $c_2(m)$ behave like $1/\sqrt m$. 
Note also that the bound in \eqref{eq:KMT} is non-trivial for $x\geq \tilde a_n:=\log(1+c_2(m)n^{1/2})/c_1(m)$.

\begin{Rem}
\label{rem:KMT}
{\rm
In order to obtain KMT bounds in our framework, we first apply Proposition~\ref{prop:KMT}, individually for each  $1\leq l\leq m$ to construct $(\widetilde{X}_{i,l}^m : 1\leq i\leq n)$, and the Brownian motion $B^{l}$ satisfying \eqref{eq:KMT} on some probability space $(\Omega_l, {\cal F}_l, \pp_l)$.
At this point, on the probability space $(\Omega_1\times\dots\times\Omega_m, {\cal F}_1\otimes\dots\otimes{\cal F}_m,\pp_1\otimes\dots\otimes\pp_m )$ the Brownian motions $B^{l}$, $1\leq l\leq m$, are rendered independent and  so are, for 
different $l$, the independent Gaussian random variables $(Y_i^l)_{1\leq i\leq n}$.
Next, we explain that this suffices and that we can consider the Brownian motion $B^{l}$, $1\leq l\leq m$, to be correlated via $\Sigma^{(m)}$, given by $\Sigma^{(m)}_{i,i}=1$, $\Sigma^{(m)}_{i,j}=-1/(m-1)$.   Indeed, setting $\overline{b}(1/m)$ for the 
centered Bernoulli distribution $b(1/m)$ normalized with unit variance, we have 
\begin{equation}
\label{eq:KMT_P}
\Big((\widetilde{X}_{i,1}^m)_{1\leq i\leq n},\dots, (\widetilde{X}_{i,m}^m)_{1\leq i\leq n}, (Y_i^1)_{1\leq i\leq n},\dots, (Y_i^m)_{1\leq i\leq n}\Big)\sim \overline{b}(1/m)^{\otimes nm}\otimes {\cal N}(0,1)^{\otimes nm}.
\end{equation}
Now, consider the vector 
\begin{equation}
\label{eq:KMT_Q}
\Big(
(U_{i,1}^m)_{1\leq i\leq n},\dots, (U_{i,m}^m)_{1\leq i\leq n}, (V_i^1)_{1\leq i\leq n},\dots, (V_i^m)_{1\leq i\leq n}\Big)\sim Q
\end{equation}
where $U_{i,l}^m\sim\overline{b}(1/m)$, $V_{i}^{l}\sim{\cal N}(0,1)$
and where both $(U_{i,l}^m)_{1\leq l\leq m}$, $(V_i^l)_{1\leq l\leq m}$ are correlated by $\Sigma^{(m)}$. 
The distribution $Q$, in \eqref{eq:KMT_Q}, is absolutely continuous with respect to the distribution given in \eqref{eq:KMT_P}, and let us denote by $f$ its 
Radon-Nikod\'ym derivative.   Then, consider the probability space $(\Omega_1\times\dots\times\Omega_m, {\cal F}_1\otimes\dots\otimes{\cal F}_m,\pp_f)$ where
$$
d{\mathbb P}_f=f\big((\widetilde{X}_{i,1}^m)_{1\leq i\leq n},\dots, (\widetilde{X}_{i,m}^m)_{1\leq i\leq n}, (Y_i^1)_{1\leq i\leq n},\dots, (Y_i^m)_{1\leq i\leq n}\big)d{\mathbb P},
$$ 
and where $\pp=\pp_1\otimes\dots\otimes \pp_m$. 
Observe that under $\pp_f$, the vector 
$$
\Big((\widetilde{X}_{i,1}^m)_{1\leq i\leq n},\dots, (\widetilde{X}_{i,m}^m)_{1\leq i\leq n}, (Y_i^1)_{1\leq i\leq n},\dots, (Y_i^m)_{1\leq i\leq n}\Big)
$$
has distribution $Q$ and, therefore, both $(\widetilde{X}_{i,l}^m)_{1\leq l\leq m}$ and 
$(Y_i^l)_{1\leq l\leq m}$ are correlated via $\Sigma^{(m)}$. 
In turn, the Brownian motions $B^l$, $1\leq l\leq m$, are also correlated via $\Sigma^{(m)}$. 
In the sequel, up to applying this argument, we assume that the Brownian motions are correlated via 
$\Sigma^{(m)}$, but we keep our notations unchanged (see also our next comment). 

Alternatively, if instead of using Proposition~\ref{prop:KMT} one uses Theorem~4.1 of \cite{L}, 
each random variable (vector) in this 
theorem should be $m$(the alphabet size)-dimensional and there should be $n$(the number of letters in the word) 
many of them.  Then, one constructs a KMT-approximation with independent 
copies of these vectors but where each vector has its coordinates taken correlated  
with the same correlation matrix as each random vector associated with the $i$th letter of the word.  More precisely, 
using 
the terminology of \cite{L}, one can construct $n$ random vectors, which are independent copies of the $X_i$'s there, 
denoted 
by $\widetilde{X}_i$, in such a way that they remain independent of each other but with the 
same covariance matrix as the $X_i$'s, as well as $n$ independent Gaussian random vectors $Y_i$'s 
with the same covariance matrix, for which a KMT 
approximation hold.   In our case, let 
$$
X_{i,j} = X_{i,j}^m - \ee[X_{i,j}^m], 
$$
where $X_{i,j}^m$ is defined at the beginning of Section 2.  Now, for fixed $i$, $X_{i,j}^m$ and $X_{i,k}^m$, $j\neq k$ are 
orthogonal random variables and so the covariance matrix $\Sigma = (\Sigma_{k,\ell})_{1\le k,l \leq m}$ of the $X_{i,j}$, where $i$ is a fixed letter of the word while $k$ and $\ell$ run over the alphabet of size $m$
is such that 
$$\Sigma_{k,\ell} = \ee[X_{i,k}X_{i,\ell}] = 
\begin{cases}
\frac{m-1}{m^2}, \quad {\rm if} \ k = \ell, \\
-\frac{1}{m^2}, \quad {\rm if} \ k\neq \ell.    
\end{cases}$$
Then, taking $X_i = (X_{i,1}, \dots, X_{i,m})$ in Theorem~4.1 of \cite{L}, one can 
further take ${\widetilde X}_i$ in the same 
theorem to have correlated entries with covariance matrix $\Sigma$.  Now, $\Sigma$
has eigenvalues $1/m$ with 
multiplicity $m$.  Therefore, the uniform non-degeneracy conditions (4.1) of Theorem~4.1 (which seems to 
contain a typo with an unnecessary extra $D^2$) are satisfied 
and so its conclusions apply.   Using this dependent version of the KMT approximation might lead to a 
different alphabet-growth rate $\alpha$, after evaluating the various parameters.  
}

\end{Rem}

\medskip
{\bf Approximating sets.}
Let $A_1^n=\{\max_{l\leq m} Y_n^{m,l}>a_n\}$, for some $a_n=Cc_1(m)^{-1}(\log n)^2\geq \tilde a_n$ 
where $C$ is some finite constant. 
We have
\begin{eqnarray*}
\pit(A_1^n)&=&\pit\Big(\bigcup_{l\leq m}\{Y_n^{m,l}>a_n\}\Big)\\
&\leq&\sum_{l\leq m}\pit(Y_n^{m,l}>a_n)\\
&\leq &m e^{-c_1(m)a_n}(1+c_2(m)n^{1/2})\\
&\sim& \sqrt{mn}e^{-c_1(m)a_n}\\
&=& \sqrt{mn}e^{-(c_3C(\log n)^2)/2}
\to 0, \quad n\to+\infty.
\end{eqnarray*}

\medskip
\noindent
Let $A_2^n=\{\max_{1\leq l\leq m} W_n^{m,l}>b_n\}$, for $b_n=\log n$. 
Standard estimates (including reflection principle, Brownian scaling and Gaussian tail estimates) lead to :
\begin{eqnarray*}
\pit(A_2^n)
&=&\pit\Big(\bigcup_{l\leq m}\{W_n^{m,l}>b_n\}\Big)\\
&\leq&\sum_{l\leq m}\pit(W_n^{m,l}>b_n)\\
&\leq&m\pit\left(W_n^{m,1}>b_n\right)\\
&=&m\pit\Big(\sup_{\begin{subarray}{c} 0\leq s,t\leq n\\ |s-t|\leq 1\end{subarray}}
|B_s^{m,1}-B_t^{m,1}|>b_n\Big).
\end{eqnarray*}
But,
\begin{eqnarray*}
\sup_{\begin{subarray}{c} 0\leq s,t\leq n\\ |s-t|\leq 1\end{subarray}}|B_s^{m,1}-B_t^{m,1}|
&\leq&\sup_{0\leq i\leq n-2}\sup_{i\leq s,t\leq i+2}|B_s^{m,1}-B_t^{m,1}|\\
&\leq&\sup_{0\leq i\leq n-2}\left(\sup_{i\leq t\leq i+2}B_t^{m,1}-\inf_{i\leq s\leq i+2}B_s^{m,1}\right),
\end{eqnarray*}
and so 
\begin{eqnarray}
\nonumber
\pit(A_2^n)
\nonumber
&\leq &m\pit\left(\sup_{0\leq i\leq n-2}\left(\sup_{i\leq t\leq i+2}B_t^{m,1}-\inf_{i\leq s\leq i+2}B_s^{m,1}\right)>b_n\right)\\
\nonumber
&\leq &mn\pit\left(\sup_{t\in [0,2]}B_t^{m,1}-\inf_{s\in[0,2]}B_s^{m,1}>b_n\right)\\
\nonumber
&\leq &mn\left(\pit\left(\sup_{t\in [0,2]}B_t^{m,1}>b_n/2\right)
+\pit\left(\sup_{s\in[0,2]}B_s^{m,1}>b_n/2\right)\right)\\
\nonumber
&\leq&2mn \pit\left(|B_2^{m,1}|>b_n/2\right)\\
\label{eq:bn2}
&\leq&4mn \exp(-b_n^2/16) \to 0,\quad n\to+\infty. 
\end{eqnarray}

\medskip
{\bf Final bound.}
Since the Brownian motions $B^l$, $1\leq l\leq m$, are correlated via $\Sigma^{(m)}$ (see Remark~\ref{rem:KMT}), combining \eqref{eq:spectrum} and \eqref{eq:GUE0} and the observation made after \eqref{eq:percolation2}, we have that $(L_k(n,m))_{1\leq k\leq r}$, for independent Brownian motions, and $(L_k(n,m))_{1\leq k\leq r}$, for Brownian motions correlated via $\Sigma^{(m)}$, only differ by $Z_m e_m$.
But since $Z_m m^{1/6}\Rightarrow 0$, the limiting result \eqref{eq:scaling2}
still applies for $(L_k(n,m))_{1\leq k\leq r}$ with Brownian motions correlated via $\Sigma^{(m)}$ (see the paragraph after \eqref{eq:nm2}). 
As a consequence, the approximation of $(\widetilde{V}_k(n,m))_{1\leq k\leq r}$ by $(L_k(n,m))_{1\leq k\leq r}$ will imply the theorem if 
\begin{equation}
\label{eq:but1}
\pit\left(\sum_{k=1}^r\left|\widetilde{V}_k(n,m)-L_k(n,m)\right|\geq c_n\right)\to 0, \quad n\to+\infty,
\end{equation}
for some 
\begin{equation}
\label{eq:cn}
c_n=o(n^{1/2}m^{-1/6}).
\end{equation}
Since $\lim_{n\to+\infty}\big(\pit(A_1^n)+\pit(A_2^n)\big)=0$, it is enough to prove that
\begin{equation}
\label{eq:aa1}
\lim_{n\to+\infty}\pit\left(\left\{\sum_{k=1}^r\left|\widetilde{V}_k(n,m)-L_k(n,m)\right|\geq c_n\right\} \cap (A_1^n)^c\cap (A_2^n)^c\right)=0.
\end{equation}
But  
\begin{align*}
&\ee\left[\sum_{k=1}^r\left|\widetilde{V}_k(n,m)-L_k(n,m)\right| \:\ind_{(A_1^n)^c\cap (A_2^n)^c}\right]\\
&\leq \sum_{k=1}^r2rm\ee\big[ (Y_n^{m,1}+W_n^{m,1}) \: \ind_{(A_1^n)^c\cap (A_2^n)^c}\big]\\
&\leq 2r^2m\left(\ee\big[ Y_n^{m,1}\: \ind_{ Y_n^{m,1}\leq a_n}\big] +b_n\right)\\
&\leq 2r^2m\left(\ee\big[ (Y_n^{m,1}-\tilde a_n)\: \ind_{Y_n^{m,1}\leq a_n}\big] +\tilde a_n+b_n\right)\\
&\leq 2r^2m\left(\ee\big[ (Y_n^{m,1}-\tilde a_n)\: \ind_{\tilde a_n \leq Y_n^{m,1}\leq a_n}\big] +\tilde a_n+b_n\right)\\
&\leq 2r^2m\left(\int_{\tilde a_n}^{a_n} \pit(Y_n^{m,1}\geq x) dx +\tilde a_n+b_n\right)\\
&\leq 2r^2m\left(\int_{\tilde a_n}^{a_n} e^{-c_1(m)x}(1+c_2(m)n^{1/2}) dx +\tilde a_n+b_n\right)\\
&\leq 2r^2m\left(\frac{1+c_2(m)n^{1/2}}{c_1(m)} e^{-c_1(m)\tilde a_n}+\tilde a_n+b_n\right)\\
&\leq 2r^2m\left(\frac 1{c_1(m)}+\tilde a_n+b_n\right)\\
&\leq 2r^2m^{3/2}\left(\frac{2(1+\log(1+c_2(m)n^{1/2}))}{c_3}+b_n\right).
\end{align*}
Finally, 
\begin{align}
\nonumber
\pit\left(\left\{\sum_{k=1}^r\left|\widetilde{V}_k(n,m)-L_k(n,m)\right|\geq c_n\right\}\cap (A_1^n)^c\cap (A_2^n)^c\right)&\\
\label{eq:aa2}
\leq \frac{2r^2m^{3/2}}{c_n}\left(\frac{2(1+\log(1+c_2(m)n^{1/2}))}{c_3}+\log n\right)=O\left(\frac{m^{3/2}\log n}{c_n}\right).
&
\end{align}
A choice of $c_n$ ensuring that the bound in \eqref{eq:aa2} goes to zero as $n\to+\infty$ and also compatible with \eqref{eq:cn}
is possible when $m^{3/2}\log n=o(n^{1/2}m^{-1/6})$, i.e., when $m=o(n^{3/10}(\log n)^{-3/5})$. 
Finally, \eqref{eq:but1} and \eqref{eq:aa1} hold true, achieving the proof of Theorem \ref{theo:main}. 
\CQFD

\begin{Rem}
\label{rem:condition}
{\rm $ $

\begin{itemize}

\item 
In the above proof, the condition $m=o(n^{3/10}(\log n)^{-3/5})$ is needed only once, to ensure the compatibility of \eqref{eq:cn} with the bound \eqref{eq:aa2}.
However, this is essential to make the Gaussian approximation work. 

\item When $m=[n^a]$, the growth condition $m=o(n^{3/10}(\log n)^{-3/5})$ rewrites as $a<3/10$, 
and this growth condition remains true, in particular, when $m$ is of sub-polynomial order. 
The condition $a<3/10$ is stronger than its counterpart $a<3/7$ in \cite{BM} 
and this seems to be due to the fact that we work with a triangular array of random variables. 

\item For the top line of the tableau, our result is  short of a result of Johansson in \cite{Johansson} which asserts 
the convergence of $V_1(n,n^a)$ (properly scaled and normalized) towards the Tracy-Widom distribution. 
More precisely, setting $a_n\ll b_n$ for $a_n=o(b_n)$, 
Th. 1.7 in \cite{Johansson} actually gives in our notations: 
for $\sqrt n\ll m$,
$$
\frac{V_1(n,m)-n/m-2\sqrt{n}}{n^{1/6}}\Rightarrow F_{TW},
$$
for $(\log n)^{3/2}\ll m \ll\sqrt n$,
$$
\frac{V_1(n,m)-n/m-2\sqrt{n}}{n^{1/2}m^{-2/3}}\Rightarrow F_{TW}, 
$$
and, for $\sqrt n/m\to l$,
$$
\frac{V_1(n,m)-n/m-2\sqrt{n}}{(1+l)^{2/3}n^{1/6}}\Rightarrow F_{TW}.  
$$
In the middle limit above, \cite[Th. 1.7]{Johansson} requires $(\log n)^{3/2}=o(m)$ 
while we do not require a lower bound condition on $m$. 
Besides, our Theorem \ref{theo:main} applies to the shape of the whole Young tableau. 
\end{itemize}
}\end{Rem}


\section{Proof of Theorem \ref{theo:nonunif_TW}}
\label{sec:nonunif_TW}

In this section, we sketch the changes needed in the previous arguments in order to prove Theorem~\ref{theo:nonunif_TW}. 
Note that in the uniform setting, the representation \eqref{eq:combinatoire_V_k} for $V_k(n,m)$ is a maximun taken over the most probable letters. 
This is trivially true since, in this case, all the letters have the same probability. 
But this property which appears to be fundamental when we center and normalize the $X_{i,j}^m$, is no longer true in the non-uniform setting. 
However, we shall approximate $V_1(n,m)$ below by a random variable $V_1'(n,m)$ defined as a maximum taken only over most probable letters as in \eqref{eq:combinatoire_V_k}, see \eqref{eq:combinatoire_V_1'}. 
Part of the remaining work is then to show that we can suitably control this approximation and this is done in Lemma \ref{lemme:control}. 
This control is at the root of the extra condition \eqref{eq:condition2} in Theorem \ref{theo:nonunif_TW}.

Let us revise our notation for the non-uniform setting. 
In this section, $X_i^m$, $1\leq i\leq n$, are  independently and identically distributed with $\pit(X_1^m=\alpha_j)=p_j^m$. 
Set $p_{max}^m=\max_{1\leq j\leq m} p_j^m$ and $J(m)=\{j\: :\: p_j^m=p_{max}^m\}=\{j_1, \dots, j_{k(m)}\}$, 
with $k(m)=\mbox{card }(J(m))$, and set also $\sigma_m^2=p_{max}^m(1-p_{max}^m)$. 
Finally, note that since $k(m(n))p_{max}^m\leq 1$ and $k(m(n))\to +\infty$, it follows that $p_{max}^{m(n)}\to 0$, as $n\to+\infty$. 
\\

{\bf Brownian scaling.}
Let $(B^{l}(s))_{s\geq 0}$, $1\leq l\leq k(m)$, be independent standard Brownian motions. 
For $s>0$, $m\geq 1$ and $k\geq 1$, let
\begin{equation}
\label{eq:percolation3}
L_1(s,k(m))=\sup_{{\bf t}\in I_{k(m)}(s)} \sum_{l=1}^{k(m)}\big(B^{l}(t_{l})-B^{l}(t_{l-1})\big),
\end{equation}
where $I_{k(m)}(s)=\{{\bf t}\: : \: 0\leq t_1\leq\dots\leq t_{l-1}\leq t_l\leq\dots\leq t_{k(m)}=s\}$. 
Recall that $L_1(1,k(m))$ has the same law as the largest eigenvalue of a $k(m)\times k(m)$ GUE matrix 
(see \eqref{eq:spectrum}, \eqref{eq:GUE0}, \eqref{eq:Johansson} and Remark \ref{rem:Johansson-TW} for $k=1$), and so: 
$$
k^{1/6}(L(1,k)-2\sqrt{k})\Rightarrow F_{TW}.
$$ 
By Brownian scaling, $L_1(s,m)\stackrel{{\cal L}}{=}\sqrt s L_1(1,m)$, so that when $n\to+\infty$: 
\begin{equation}
\label{eq:scaling3}
\frac{L_1(n,k(m(n)))-2\sqrt{nk(m(n))}}{n^{1/2}k(m(n))^{-1/6}}
\Rightarrow F_{TW}.
\end{equation}

{\bf Combinatorics revisited.}
Let 
$$
X_{i,j}^m=\left\{\begin{array}{ll}
1&\mbox{when } X_i^m=\alpha_j\\
0&\mbox{otherwise,}
\end{array} \right .
$$
be Bernoulli random variables with parameter $\pit(X_i^m=\alpha_j)=p^m_j$ and variance $(\sigma_j^m)^2=p_j^m(1-p_j^m)$.
For a fixed $1\leq j\leq m$, the $X_{i,j}^m$s are independent and identically distributed. 
Since the expression \eqref{eq:combinatoire_V_1} has a purely combinatorial nature, we still have 
\begin{eqnarray*}
V_1(n,m)
&=&\max_{\begin{subarray}{l}0=l_0\leq l_1\leq\cdots\\\leq l_{m-1}\leq l_m=n\end{subarray}}
\Big(\sum_{j=1}^m \sum_{i=l_{j-1}+1}^{l_j}X_{i,j}^m\Big),
\end{eqnarray*}
with the convention that $\sum_{i=l_{j-1}+1}^{l_j}X_{i,j}^m=0$, whenever $l_{j-1}=l_j$. 

In fact, for most draws the maximum in $V_1$ is attained on the sums $\sum_{j\in J(m) } \sum_{i=l_{j-1}+1}^{l_j}X_{i,j}^m$ corresponding to the most probable letters, 
that is, letting 
\begin{equation}
\label{eq:combinatoire_V_1'}
V_1'(n,m)
=\max_{\begin{subarray}{l}0=l_0\leq l_1\leq\cdots\\\leq l_{m-1}\leq l_m=n\\l_{j-1}=l_j \ {\rm for } \ j\not\in J(m)\end{subarray}}
\Big(\sum_{j=1}^n \sum_{i=l_{j-1}+1}^{l_j}X_{i,j}^m\Big),
\end{equation}
we have, with large probability, $V_1(n,m)=V_1'(n,m)$. 
However, it is not always true that $V_1(n,m)=V_1'(n,m)$, 
for instance if the $n$ letters drawn are letters with associated 
probability strictly less than $p_{max}^m$, $V_1'(n,m)=0$ 
while there is a $l=(l_j)_{j=0,\dots, m}$ with $0=l_0\leq l_1\leq\cdots\leq l_{m-1}\leq l_m=n$ such that $\sum_{j=1}^m \sum_{i=l_{j-1}+1}^{l_j}X_{i,j}^m>0$, 
ensuring that $V_1(n,m)>0$.  
In the sequel, we prove Theorem \ref{theo:nonunif_TW} by first showing that the statement of the theorem is true for $V_1'(n,m)$ instead of $V_1(n,m)$ 
and then by controling the error made when $V_1'(n,m)$ is replaced by $V_1(n,m)$. 

\bigskip

{\bf Centering and reducing.}
Let $\widetilde{X}_{i,j}^m=(X_{i,j}^m-p_j^m)/\sigma_j^m$ be the corresponding centered and normalized scaled Bernoulli random variables
and let $\widetilde{S}_l^{m,j}=\sum_{i=1}^l \widetilde{X}_{i,j}^{m}$. 
Let also,
\begin{eqnarray*}
\widetilde V_1'(n,m)
&=&\max_{\begin{subarray}{l}0=l_0\leq l_{j_1}\leq\cdots\\\leq l_{j_{k(m)-1}}\leq l_{j_{k(m)}}=n\end{subarray}}\left(\sum_{j\in J(m)} \sum_{i=l_{j-1}+1}^{l_j}\widetilde X_{i,j}^m\right)\\
&=&\max_{\begin{subarray}{l}0=l_0\leq l_{j_1}\leq\cdots\\\leq l_{j_{k(m)-1}}\leq l_{j_{k(m)}}=n\end{subarray}}
\left(\sum_{j\in J(m)} (\widetilde S_{l_j}^{m,j}-\widetilde S_{l_{j-1}}^{m,j})\right)\\
&=&\sup_{{\bf t}\in I_{k(m(n))}(n)} \sum_{\delta=1}^{k(m(n))}\left(\widetilde S_{[t_{j,l}]}^{m,j_\delta}-\widetilde S_{[t_{j,l-1}]}^{m,j_\delta}\right),
\end{eqnarray*}
which is to be compared to \eqref{eq:percolation3}. 
Since $V_1'(n,m)-np_{max}^m=\sigma_m \widetilde V_1'(n,m)$, we have 
\begin{align*}
&k(m)^{1/6}\frac{V_1'(n,m)-np_{max}^m-2\sqrt{nk(m)\sigma_m^2}}{\sqrt{p_{max}^mn}}\frac{\sqrt{p_{max}^m}}{\sigma_m}\\
&\hskip 4cm = k(m)^{1/6}\frac{\widetilde V_1(n,m)-2\sqrt{nk(m)}}{\sqrt n}.
\end{align*}
Since $\sigma_m\sim \sqrt{p_{max}^m}$ and 
\begin{eqnarray*}
\frac{2\sqrt{k(m)np_{max}^m}-2\sqrt{k(m)n\sigma_m^2}}{\sqrt{n\sigma_m^2}}
&=&\frac{2\sqrt{k(m)}}{\sigma_m} \frac{p_{max}^m-\sigma_m^2}{\sqrt{p_{max}^m}+\sqrt{\sigma_m^2}}\\
&\sim&\frac{2\sqrt{k(m)}}{\sqrt{p_{max}^m}} \frac{(p_{max}^m)^2}{\sqrt{p_{max}^m}}\\
&\leq& 2\sqrt{p_{max}^m}\to 0, \quad n\to+\infty,
\end{eqnarray*}
it remains to show that
\begin{equation}
\label{eq:but11}
k(m)^{1/6}\frac{\widetilde V_1'(n,m)-2\sqrt{nk(m)}}{\sqrt n}\Rightarrow F_{TW}
\end{equation}
for which we shall use \eqref{eq:scaling3}. 

\bigskip

{\bf Sketch of proof of \eqref{eq:but11}.}
Roughly speaking, the proof of \eqref{eq:but11} is along the same lines of the corresponding proof of the convergence of \eqref{eq:V_k}, changing only $m$ into $k(m)$. 
We show that when $k(m(n))=o(n^{3/10}(\log n)^{-3/5})$, then for some Brownian motions given via strong approximation, 
we have 
$$
\left|\widetilde{V}_1'(n,m)-L_1(n,k(m(n)))\right|\\
\leq \sum_{l=1}^{k(m(n))}\left(Y_n^{m,l}+W_n^{m,l}\right),
$$
where
$$
Y_n^{m,l}=\max_{1\leq i\leq n} |S_i^{m,l}-B^{m,l}(i)|
\quad \mbox{ and } \quad 
W_n^{m,l}=\sup_{\begin{subarray}{c} 0\leq s,t\leq n\\ |s-t|\leq 1\end{subarray}}|B^{m,l}(s)-B^{m,l}(t)|.
$$

Indeed, setting $A_1^n=\{\max_{l\leq k(m(n))} Y_n^{m,l}>a_n\}$, for some $a_n=O(c_1(k(m(n)))^{-1}(\log n)^2)\geq \tilde a_n:=\log(1+c_2(k(m(n)))n^{1/2})/c_1(k(m(n)))$,
and setting $A_2^n=\{\max_{1\leq l\leq k(m(n))} W_n^{m,l}>b_n\}$, for some $b_n=O(\log n)$, 
we show that
$$
\pit(A_1^n)\to 0, \quad \pit(A_2^n)\to 0, \quad \mbox{ when } n\to +\infty.
$$ 
From \eqref{eq:scaling3}, the approximation of $\widetilde{V}_1(n,k(m(n)))$ 
by $L_1(n,k(m(n))))$ will imply the theorem if 
\begin{equation}
\label{eq:but1b}
\pit\left(\left|\widetilde{V}'_1(n,k(m(n)))-L_1(n,k(m(n)))\right|\geq c_n\right)\to 0, \quad n\to+\infty,
\end{equation}
for some 
\begin{equation}
\label{eq:cn2}
c_n=o(n^{1/2}k(m(n))^{-1/6}).
\end{equation}
Since $\lim_{n\to+\infty}\big(\pit(A_1^n)+\pit(A_2^n)\big)=0$
and
\begin{align}
\nonumber
\pit\left(\left\{\left|\widetilde{V}'_1(n,k(m(n)))-L_1(n,k(m(n)))\right|\geq c_n\right\} \cap (A_1^n)^c\cap (A_2^n)^c\right)&\\
\label{eq:cn3}
\leq \frac{2k(m(n))^{3/2}}{c_n}\left(\frac{2(1+\log(1+c_2(k(m(n)))n^{1/2}))}{c_3}+\log n\right),&
\end{align}
a choice of $c_n$, ensuring that the bound in \eqref{eq:cn3} goes to zero 
and is compatible with \eqref{eq:cn2}, is possible since $k(m(n))=o(n^{3/10}(\log n)^{-3/5})$.
This proves \eqref{eq:but11} and thus the statement \eqref{eq:nonunif_TW} of 
Theorem \ref{theo:nonunif_TW}, but for $V_1'(n,m)$ instead of $V_1(n,m)$.

\bigskip

{\bf Control of the error $V_1(n,m)-V_1'(n,m)$.}
Clearly $V_1(n,m)-V_1'(n,m)\geq 0$ and is, in fact, zero with a large probability, so that we expect $\ee[V_1(n,m)-V_1'(n,m)]$ to be small. 
Actually we show:
\begin{lemme}
\label{lemme:control}
For some absolute constant $C>0$, we have 
\begin{equation}
\label{eq:control}
\ee[|V_1(n,m)-V_1'(n,m)|]\leq Cnp_{2nd}^m,
\end{equation}
where $p_{2nd}^m$ stands for the second largest probability for the letters of ${\cal A}_m$.
\end{lemme}
The conclusion in \eqref{eq:nonunif_TW} holds true when 
\begin{equation}
\label{eq:control2}
\lim_{n\to+\infty}\left(\ee[|V_1(n,m)-V_1'(n,m)|] \times 
\frac{k(m(n))^{2/3}}{\sqrt{k(m(n))p_{max}^{m(n)}n}}\right)=0.
\end{equation}
But with the help of \eqref{eq:control}, the conclusion 
in \eqref{eq:control2} is then valid 
when $\lim_{n\to+\infty}\frac{p_{2nd}^{m(n)}k(m(n))^{1/6} n^{1/2}}{(p_{max}^{m(n)})^{1/2}}=0$ 
and, since $k(m(n))=o(n^{3/10}(\log n)^{-3/5})$, this will 
follow from \eqref{eq:condition2}. 

It remains to prove Lemma \ref{lemme:control}, i.e. to give an explicit bound on $\ee[|V_1(n,m)-V_1'(n,m)|]$.
To do so, rewrite $V_1(n,m)=\max_{l\in I(m)} Z(l)$ and $V_1'(n,m)=\max_{l\in I^*(m)} Z(l)$ where
$I(m)=\{l=(l_j)_{1\leq j\leq m}\: : \: l_{j-1}\leq l_j, l_0=0, l_m=n\}$, 
$I^*(m)=\{l\in I(m) : l_{j-1}=l_j \mbox{ for } j\not\in J(m)\}$
and 
$$
Z(l)=\sum_{j=1}^m Y_j(l), \quad Y_j(l)=\sum_{i=l_{j-1}+1}^{l_j} X_{i,j}^m.
$$
Clearly, since $I^*(m)\subset I(m)$, we have $V_1'(n,m)\leq V_1(n,m)$. 
Moreover, since the $X^m_{i,j}$ are Bernoulli random variables with parameter $p_j^m$ 
and since the $X_i$s are independent, we have $Y_j(l)\sim {\cal B}(l_j-l_{j-1}, p_j)$ 
and $\sum_{j\in J(m)} Y_j(l)\sim {\cal B}(\sum_{j\in J(m)}l_j-l_{j-1}, p_{max}^m)$, where ${\cal B}(n,p)$ stands for the binomial distribution with parameters $n$ and $p$. 

\bigskip\noindent
If $l\in I^*(m)$, $Z(l)=\sum_{j\in J(m)} Y_j(l)\sim {\cal B}(n, p_{max}^m)$ since in this case
$n=\sum_{j=1}^m (l_j-l_{j-1})=\sum_{j\in J(m)}(l_j-l_{j-1})$. 
If $l\not\in I^*(m)$, we rewrite $Z(l)$ as:
$$
Z(l)=Z(\tilde l)+ R(l),
$$
where $\tilde l\in I^*(m)$ and $R(l)$ is an error term. 
Indeed, let $J_l=\{j\not\in J(m)\: :\: l_{j-1}<l_j\}$ and for $j\in J_l$, define
$$
\theta(j)=\left\{
\begin{array}{ll}
\max A_j &\mbox{ if } A_j\neq \emptyset\\
\min B_j &\mbox{ otherwise,}
\end{array}
\right .
$$
where $A_j=\{k\in J(m)\: : \: k\leq j\}$ and where $B_j=\{k\in J(m)\: : \: k\geq j\}$. 
Now,
\begin{eqnarray}
\nonumber
Z(l)&=&
\sum_{j\in J(m)}\sum_{i=l_{j-1}+1}^{l_j} X_{i,j}^m
+\sum_{j\in J_l}\sum_{i=l_{j-1}+1}^{l_j} X_{i,j}^m\\
\label{eq:Ztilde}
&=&
\sum_{j\in J(m)}\sum_{i=l_{j-1}+1}^{l_j} X_{i,j}^m
+\sum_{j\in J_l}\sum_{i=l_{j-1}+1}^{l_j} X_{i,\theta(j)}^m\\
\label{eq:R}
&&+\sum_{j\in J_l}\sum_{i=l_{j-1}+1}^{l_j} (X_{i,j}^m-X_{i,\theta(j)}^m). 
\end{eqnarray}
Define $\tilde l\in I^*(m)$ by $\tilde l_j=\tilde l_{j-1}$ if $j\not\in J(m)$ and 
$\tilde l_j=l_{k-1}$ for $j\in J(m)$, where $k=\min\{ l>j\: : \: l\in J(m)\}$, with the convention that $\min\emptyset=m+1$, 
and that $\tilde l_{j_0-1}=0$, for $j_0=\min J(m)$. 
We then have 
$$
\sum_{j\in J(m)}\sum_{i=l_{j-1}+1}^{l_j} X_{i,j}^m
+\sum_{j\in J_l}\sum_{i=l_{j-1}+1}^{l_j} X_{i,\theta(j)}^m
=Z(\tilde l).
$$
\medskip

\noindent 
Let $\alpha_{i,j}^m:=X_{i,j}^m-X_{i,\theta(j)}^m$ be the random variables taking the values $-1, 0$ and $+1$ with respective probabilities 
$p_{max}^m, 1-p_{max}^m-p_j^m$ and $p_j^m$. 
Independently, let $\epsilon_{i,j}^m$ be Bernoulli random variables with parameter $q_j^m=(p_{2nd}^m-p_j^m)/(1-p_{max}^m-p_j^m) \in (0,1)$ 
where $p_{2nd}^m=\max(p_j^m<p_{max}^m\: : \: 1\leq j\leq m)$ and define
$$
\beta_{i,j}^m=
\left\{
\begin{array}{rl}
-1 & \alpha_{i,j}^m=-1\\
0 & \alpha_{i,j}^m=0 \mbox{ and } \epsilon_{i,j}^m=0\\
+1 & \alpha_{i,j}^m=+1 \mbox{ or } \alpha_{i,j}^m=0 \mbox{ and } \epsilon_{i,j}^m=1.
\end{array}
\right .
$$
Note that $\pit(\beta_{i,j}^m=+1)=p_{2nd}^m$ and that $\alpha_{i,j}^m\leq \beta_{i,j}^m$,
so that 
$$
R(l)\leq  \widetilde R(l)=\sum_{j\in J_l}\sum_{i=l_{j-1}+1}^{l_j} \beta_{i,j}^m. 
$$
Since $Z(l)\leq Z(\tilde l)+\widetilde R(l)$, we have 
\begin{eqnarray*}
\max_{l\in I(m)} Z(l)&\leq& \max_{l\in I(m)} Z(\tilde l)+\max_{l\in I(m)} \widetilde R(l)\\
&\leq& \max_{l\in I^*(m)} Z(l)+\max_{l\in I(m)} \widetilde R(l).
\end{eqnarray*}
Next, observe that for $l\in I^*(m)$, $\widetilde R(l)=0$. 
However since the event $\{\widetilde R(l)<0, \: \forall l\not\in I^*(m)\}$ 
is non-negligible, we cannot change $\max_{l\in I(m)} \widetilde R(l)$ 
into $\max_{l\not\in I^*(m)} \widetilde R(l)$. 
We obtain
$$
0\leq \max_{l\in I(m)} Z(l)- \max_{l\in I^*(m)} Z(l)\leq \max_{l\in I(m)} \widetilde R(l).
$$
The random variable $\widetilde R(l)$ is the sum of $\sum_{j\in J_l} (l_j-l_{j-1})$ {\it iid} 
random variables, 
so that $\max_{l\in I(m)} \widetilde R(l)$ is distributed according to 
$\left(\max_{1\leq k\leq n} \sum_{i=1}^k \beta_i^m\right)^+$
where $(\beta_i^m)_i$ are {\it iid} with 
\begin{equation}
\label{eq:proba_ponctuelle}
\pit(\beta_1^m=-1)=p_{max}^m, \quad 
\pit(\beta_1^m=0)=1-p_{max}^m-p_{2nd}^m,\quad 
\pit(\beta_1^m=+1)=p_{2nd}^m. 
\end{equation}
\medskip 
\noindent
We are now interested in bounding $\ee\left[\left(\max_{1\leq k\leq n} 
\sum_{i=1}^k \beta_i^m\right)^+\right]$. 

\noindent
Let $(\epsilon_i^m)_i$ be {\it iid} Bernoulli random variables with parameter $p_{max}^m+p_{2nd}^m$ 
and let, independently,  $(Y_i^m)_i$ be {\it iid} Rademacher random variables with parameter 
$p_{2nd}^m/(p_{2nd}^m+p_{max}^m)$ (i.e. $\pit(Y_i^m=1)=1-\pit(Y_i^m=-1)=p_{2nd}^m/(p_{2nd}^m+p_{max}^m)$).
Then $\beta_i^m$ and $\epsilon_i^m Y_i^m$ have the same distribution
and we have
$$
\ee\left[\left(\max_{1\leq k\leq n}\sum_{i=1}^k\beta_i^m\right)^+\right]
=\ee\left[\left(\max_{1\leq k\leq n}\sum_{i=1}^k\epsilon_i^m Y_i^m\right)^+\right]
=\ee\left[\ee\left[\left(\max_{1\leq k\leq n}\sum_{i=1}^k\epsilon_i^m Y_i^m\right)^+\Big|{\cal G}_n\right]\right]
$$
where ${\cal G}_n=\sigma(\epsilon_i^m\: :\: 1\leq i\leq n)$. 
But since $(\epsilon_i^m)_i$ is independent of $(Y_i^m)_i$, we have 
$$
\ee\left[\left(\max_{1\leq k\leq n}\sum_{i=1}^k\epsilon_i^m Y_i^m\right)^+\Big|{\cal G}_n\right]
=\ee\left[\left(\max_{1\leq k\leq \ell}\sum_{i=1}^k Y_i^m\right)^+\Big|{\cal G}_n\right]
$$
where $\ell=\sum_{i=1}^n \epsilon_i^m$ has a ${\cal B}(n,p_{max}^m+p_{2nd}^m)$ distribution. 
But 
\begin{eqnarray*}
\ee\left[\left(\max_{1\leq k\leq \ell} \sum_{i=1}^k Y_i^m\right)^+\Big|{\cal G}_n\right]
&=&\sum_{k=1}^{+\infty} \left(\pit\left(\left(\max_{1\leq j\leq \ell} \sum_{i=1}^j Y_i^m\right)^+\geq k\right)\right)\\
&=&\sum_{k=0}^{+\infty}\Big( 1-\pit\big(\max_{1\leq j\leq \ell} \sum_{i=1}^j Y_i^m\leq k\big)\Big)\\
&=&\sum_{k=0}^{\ell-1} \Big(1-\pit\big(\max_{1\leq j\leq \ell} \sum_{i=1}^j Y_i^m\leq k\big)\Big)\\
&=&\ell-U_\ell
\end{eqnarray*} 
where $U_\ell=\sum_{k=0}^{\ell-1} u_{\ell,k}$ and 
$u_{\ell,k}=\pit(\max_{1\leq j\leq \ell} \sum_{i=1}^j Y_i^m\leq k)$. 
With the latest notation, we are now investigating $\gamma_n=\ee[\ell-U_\ell]$. 
For simplicity, in the sequel, we set $p_{*,m}:=p_{2nd}^m/(p_{2nd}^m+p_{max}^m)$ and $q_{*,m}:=1-p_{*,m}$. 

The elements of the sequence $(u_{\ell,k})_{1\leq k\leq \ell-1}$ satisfy the following induction relations:
\begin{align*}
&u_{\ell,k}=q_{*,m}u_{\ell-1,k+1}+p_{*,m}u_{\ell-1, k-1},\: k\geq 1, \quad u_{\ell,0}=q_{*,m}u_{\ell-1,1},
\end{align*}
and $u_{\ell,k}=1$ for $k\geq \ell$. 
From it, we derive $U_\ell=2q_{*,m}-q_{*,m}u_{\ell-1,0}+U_{\ell-1}$
and, since $U_1=u_{1,0}=q_{*,m}$, 
$U_\ell=(2\ell-1)q_{*,m}-q_{*,m}\sum_{k=1}^{\ell-1} u_{k,0}$.

In order to compute $\sum_{k=1}^{\ell-1} u_{k,0}$, we introduce the hitting time 
$\tau_1^m=\min\left(k\geq 1\: : \: \sum_{i=1}^k Y_i^m=1\right)$ of the random walk $(\sum_{i\leq j} Y_i^m)_j$. 
We then have
$$
\pit(\tau_1^m\leq k)=\pit\Big(\max_{i\leq k}\sum_{j=1}^i Y_j^m\geq 1\Big)
=1-\pit\Big(\max_{i\leq k} \sum_{j=1}^i Y_j^m\leq 0\Big)
=1-u_{k,0},
$$
so that 
$\sum_{k=1}^{\ell-1} u_{k,0}=\sum_{k=1}^{\ell-1} \pit(\tau_1^m\geq k+1)=\sum_{k=2}^{\ell} \pit(\tau_1^m\geq k)
=-1+\sum_{k=1}^\ell \pit(\tau_1^m\geq k)$ and 
\begin{eqnarray*}
U_\ell&=&2\ell q_{*,m}-q_{*,m}\sum_{k=1}^\ell \pit(\tau_1^m\geq k)\\
&=&2\ell q_{*,m}-q_{*,m}\sum_{i=1}^{+\infty}(i\wedge\ell) \pit(\tau_1^m=i)\\
&=&2\ell q_{*,m}-q_{*,m}\ee[\tau_1^m\wedge\ell |{\cal G}_n].
\end{eqnarray*}
Next, 
$$
\ee\left[\left(\max_{i\leq k\leq \ell}\sum_{i=1}^k Y_i^m\right)^+\Big|{\cal G}_n\right]
=\ell(1-2q_{*,m})+q_{*,m}\ee[\tau_1^m\wedge\ell|{\cal G}_n],
$$
and we have 
\begin{eqnarray}
\nonumber
\gamma_n&:=&\ee[\ell-U_\ell]\\
&=&\ee\big[\ell(1-2q_{*,m})+q_{*,m}\ee[\tau_1^m\wedge\ell |{\cal G}_n]\big]\\
\nonumber
&=&\ee[\ell(1-2q_{*,m})]+q_{*,m}\ee[\tau_1^m\wedge\ell]\\
\nonumber
&=&\ee[\ell(1-2q_{*,m})]+q_{*,m}(\ee[\ell : \tau_1^m\geq n]+\ee[\tau_1^m\wedge\ell :\tau_1^m<n])\\
\label{eq:alpha3}
&=&\ee[\ell(1-2q_{*,m})]+q_{*,m}(\ee[\ell]\pit(\tau_1^m\geq n)+\ee[\tau_1^m\wedge\ell :\tau_1^m<n]).
\end{eqnarray}
From \cite[p. 352]{Feller}, the law of $\tau_1^m$ is given for $j\geq 0$ by $\pit(\tau_1^m=2j)= 0$ and
$$
\pit(\tau_1^m=2j+1)= \frac 1{2j+1} {2j+1\choose j+1} p_{*,m}^{j+1}q_{*,m}^j.
$$
Observe that the sum 
$\sum_{j=0}^{+\infty}\frac 1{2j+1}{2j+1\choose j+1} p_{*,m}^jq_{*,m}^j$
converges uniformly with respect to $n$ since $q_{*,m}\leq 1$ and $p_{*,m}\searrow 0$, when $m=m(n)\to+\infty$
so that,
\begin{eqnarray*}
\lim_{n\to+\infty}
\pit(\tau_1^{m(n)}<n)&=&\lim_{n\to+\infty}\left(\Big(\sum_{j=0}^{n-1}
\frac 1{2j+1}{2j+1\choose j+1}p_{*,m}^jq_{*,m}^j\Big) p_{m,*}\right)
=0,
\end{eqnarray*}
and $\pit(\tau_1^m\geq n)= 1$. 
Taking $n\to+\infty$ in \eqref{eq:alpha3} and recalling the definition of $p_{m,*}$, $q_{m,*}$, 
we obtain $\gamma_n\sim np_{2nd}^m$, which concludes the proof of Lemma \ref{lemme:control}. 
\CQFD


\section{Concluding remarks}
\label{sec:remarks} 

A natural question to handle next would be to derive a result similar to Theorem \ref{theo:main} for non uniformly distributed letters. 
The special case of the longest increasing subsequence (i.e., $r=1$) is dealt with in Theorem~\ref{theo:nonunif_TW}. 
Let us investigate what happens for the whole shape of the Young tableau. 

First, let us slightly expand our notation. 
In this section, $X_i^m$, $1\leq i\leq n$, are 
independently and identically distributed with $\pit(X_1^m=\alpha_j)=p_j^m$. 
In order to simplify the notations, we assume (without loss of generality) 
that the ordered letters $\alpha_1^m<\dots<\alpha_m^m$ have moreover non-increasing probabilities 
 (i.e. $p_1^m\geq p_2^m\geq \dots\geq p_m^m$). 
Let $d_i^m=\mbox{card }\{ j\::\: p_j^m=p_i^m\}$ be the multiplicity of $p_i^m$
and let $m_r^m=\max\{i\: :\: p_i^m >p_r^m\}$ be the number of letters 
(strictly) more probable than $\alpha_r^m$.
Let $J_r(m)=\{i\::\: p_i^m=p_r^m\}=\{m_r+1, \dots, m_r+d_r^m\}$ be the indices of the letters with probability $p_r^m$. 
We recover our previous notation, $r=1$, with $k(m)=d_1^m$ and $J(m)=J_1(m)$. 
Since the expression \eqref{eq:combinatoire_V_k} has a purely combinatorial nature, it still holds true that 
\begin{eqnarray*}
V_r(n,m)=\max_{{\bf k}\in J_{r,m}(n)}\left(\sum_{j=1}^r\sum_{l=j}^{m-r+j}\sum_{i=k_{j,l-1}+1}^{k_{j,l}} X_{i,l}^m\right).
\end{eqnarray*}
Let $\nu_k^m=\sum_{i=1}^kp_i^m$. 
Note that from Theorem 5.2 in \cite{HL2}, when $m$ is fixed and $ n\to+\infty$, we have 
for each $1\leq r\leq m$:
\begin{equation}
\label{eq:HL2bis}
\left(\frac{V_k(n,m)-\nu_k^mn}{\sqrt n}\right)_{1\leq k\leq r} \Longrightarrow (V_\infty^k)_{1\leq k\leq r},
\end{equation}
where the limit is given in Section 6 of \cite{HL2} by $V_\infty^r=Z(m,r)+\sqrt{p_r^m}D_{r-m_r^m,d_r^m}$,
with $Z(m,r)\sim {\cal N}(0,v_r^m)$, 
for $v_r^m=\nu_{m_r^m}^m(1-\nu_{m_r^m}^m)+(p_r^m(r-m_r^m))^2$ and, 
$$
D_{r,m}=\max_{{\bf t}\in I_{r,m}}\left( \sum_{j=1}^r\sum_{l=j}^{(m-r+j)} \big(B^l(t_{j,l})-B^l(t_{j,l-1})\big)\right),
$$
for 
\begin{eqnarray*}
I_{r,m}&=&\big\{ {\bf t}=(t_{j,l}, 1\leq j\leq r,\ 0\leq l\leq m) \: :\: 
t_{j,j-1}=0, t_{j,m-r+j}=1, 1\leq j\leq r,\\
&& t_{j,l-1}\leq t_{j,l}, 1\leq j\leq r, 1\leq l\leq m-1, t_{j,l}\leq t_{j-1, l}, 
2\leq j\leq r, 1\le l\leq m-1\big\}.
\end{eqnarray*}
Note that $D_{r,m}$ is a natural generalization of the Brownian 
functional $L_1(s,k)$ used in Section~\ref{sec:nonunif_TW}
(see also, in a queuing context, \cite{GW} and \cite{Ba}). 
In particular, $D_{r,m}$ is equal in distribution 
to the sum of the $r$ largest eigenvalues of an $m\times m$ matrix from the GUE
and Theorem~\ref{theo:Johansson} rewrites as 
\begin{equation}
\label{eq:Dkm}
\big(m^{1/6}(D_{k,m}-k\sqrt{m})\big)_{1\leq k\leq r} \Rightarrow {\bf F}_{r}{\bf\Theta}_r^{-1}, 
\quad m\to+\infty.
\end{equation}
Arguing like in the previous sections, we would like to derive the 
fluctuations of $(V_k(n,m))_{1\leq k\leq r}$ with respect to $n$ and $m$ 
simultaneously from \eqref{eq:HL2bis} and \eqref{eq:Dkm}.
However, in the non-uniform case, this is not that transparent since, for each $r\geq 1$, 
the behavior of $m_r^m$ and of $d_r^m$, with respect to $m$, is not that clear cut .
In particular, $r-m_r^m$ may not be stationnary and $\eqref{eq:Dkm}$ can no longer be 
used for $D_{r-m_r^m, d_r^m}$. 
Besides, the random fluctuations of $\sqrt{p_r^m}D_{r-m_r^m,d_r^m}$ in $V_\infty^r$ are 
of order $(p_r^m)^{1/2}(d_r^m)^{1/6}$ 
which, in general, does not dominate those of $Z(m,r)\sim {\cal N}(0,v_r^m)$. 
Thus, for general non-uniform alphabets, we cannot infer which 
part of the law of $V_\infty^r=Z(m,r)+\sqrt{p_r^m}D_{r-m_r^m,d_r^m}$ will drive the fluctuations. 
We can imagine that, taking simultaneous limits in $n$ and $m$, 
the fluctuations of $V_r(n,m(n))$, properly centered and normalized, are either 
Gaussian, either driven by ${\bf F}_r$ as in Theorem~\ref{theo:main}, 
or given by an interpolation between these distributions, depending on the alphabets considered. 


\bigskip\noindent{\bf Acknowledgments.} 
The first author thanks the School of Mathematics of the Georgia Institute 
of Technology where this work was initiated. 
The second author would like to thank the \'Equipe Modal'X of the University Paris X for its hospitality and support while part of this research was performed.  Both authors would like to Boris Bukh and Alperen {\H O}zdemir 
for their comments which led to the extended Remark 8 of the present version.  


\end{document}